\documentclass[a4paper, 11pt]{article} 
\usepackage{amsfonts,amsmath,amssymb,amscd}
\usepackage{latexsym}
\newcommand{\disp}{\displaystyle}

\newcommand{\bi}{\mathcal I}
\newcommand{\bo}{\mathcal O}
\newcommand{\bz}{\mathcal Z}
\newcommand{\bj}{\mathcal J}
\newcommand{\bh}{\mathcal H}

\newcommand{\bm}{\mathcal{M}}
\newcommand{\ba}{\mathcal {A}}
\newcommand{\bt}{\mathcal {T}}
\newcommand{\bb}{\mathcal {B}}

\newcommand{\Nn}{\mathbb{N}}
\newcommand{\Rr}{\mathbb{R}}

\newcommand{\Zz}{\mathbb{Z}}
\newcommand{\Cc}{\mathbb{C}}

\newcommand{\Qq}{\mathbb{Q}}

\newcommand{\conj}{\text{Conj}}

\newcommand{\si}{\sigma}

\newcommand{\kah}{\text{K\"{a}hler }}
\newcommand{\del}{\bar \partial}

\newcommand{\Rhd}{\Rr H^0(X,L^d)}
\newcommand{\hd}{H^0(X,L^d)}
\newcommand{\csi}{C_{\sigma}}
\newcommand{\Rcsi}{\Rr C_{\sigma}}
\newcommand{\RX}{X \setminus \Rr X}

\newcommand{\beq}{\begin{eqnarray*}}
\newcommand{\eeq}{\end{eqnarray*}}
\newcommand{\bpr}{\begin{preuve}}
\newcommand{\epr}{\end{preuve}}

\newenvironment{preuve}[1][]
{\vskip 2mm  {\it \bf Proof#1. }}{$\Box$ \vskip 2mm}

\newtheorem{remark}{Remark}
\newtheorem{defi}{Definition}
\newtheorem{theo}{Theorem}
\newtheorem{lemme}{Lemma}
\newtheorem{prop}{Proposition}
\newtheorem{coro}{Corollary}

\title{Exponential rarefaction\\ of  real  curves with many components}
\author{Damien Gayet, Jean-Yves Welschinger}

\begin{document}
\large
\maketitle
\centerline{\textbf{Abstract}}
Given a positive real Hermitian holomorphic line bundle $L$ over a smooth real projective manifold $X$,
the space of real holomorphic sections of the bundle $L^d$ inherits for every $d\in \Nn^*$
a $L^2$ scalar product which induces a Gaussian measure. When $X$ is a curve or a surface, 
we estimate the volume
of the cone of real sections whose vanishing locus contains many real components. 
In particular, the volume of the cone of maximal real  sections decreases exponentially
as $d$ grows to infinity.\\


\section*{Introduction}

Let $(X,c_X)$ be a smooth real projective manifold of dimension $n$
and $(L,c_L) \overset{\pi}{\rightarrow} (X,c_X)$ be a real ample holomorphic line bundle. In particular, $c_X$ and $c_L$ are antiholomorphic involutions on $X$ and $L$ respectively,
such that $c_X \circ \pi = \pi \circ c_L $. Let $h$ be a real Hermitian metric
on $(L,c_L)$ with positive curvature $\omega$. It induces a K\"ahler structure
on $(X,c_X)$.  For every nonnegative integer $d$,
this metric induces a Hermitian metric $h^d$ on $L^d$
and then a $L^2$-Hermitian product on the complex
vector space $H^0(X,L^d)$ of holomorphic sections of $L^d$.
This product is defined by 
$(\si, \tau) \in  H^0(X,L^d) \times H^0(X,L^d) \mapsto \langle \si,\tau \rangle = \int_X h^d(\sigma, \tau) dx \in \Cc$,
where $dx = \omega^n/\int_X \omega^n$
is the normalized volume induced by the $\kah$ form. Let $\Rr H^0(X,L^d)$ be the space of real sections $\{\si \in  H^0(X,L^d) \ | \ c_L\circ \si = \si \circ c_X \}$
and 
$\Delta_k \subset H^0(X,L^d)$ (resp.
$\Rr \Delta_k \subset \Rr H^0(X,L^d)$) be the discriminant locus (resp. its real part), that is
the set of sections which do not vanish transversally. 
For every $\si \in   H^0(X,L^d)\setminus \{0\}$, denote by $C_\sigma = \sigma^{-1}(0)$
the vanishing locus of $\sigma$ and when $\sigma$
is real, by $\Rr C_{\sigma}$ its real part.
The divisor $C_\si$ is smooth whenever $\sigma \in  \Rr H^0(X,L^d) \setminus \Rr \Delta_d$. In this case, 
we denote by $b_0(\sigma) = b_0(\Rr C_{\sigma})$ 
the number of connected components of $\Rr C_{\sigma}$.

\subsection{Real projective surfaces} 
When  $X$ is two-dimensional, we know from  Harnack-Klein 
inequality \cite{Harnack}, \cite{Klein} that
$  b_0(\Rr C_{\sigma}) \leq g(C_{\sigma}) +1,$ where
 equality holds for the so-called maximal curves.
Here, the genus $g(\csi)$ of these smooth curves $\csi$
gets computed by the adjunction formula and 
equals  $ g(\csi) = \frac{1}{2}(d^2 L^2 - d c_1(X). L +2),$
where $c_1(X)$ denotes the first Chern class of the surface $X$. 
For every $d\in \Nn^*$ and $a \in \Qq_+^*$, denote by
$$ \bm_d^a = \{ \si \in \Rhd \setminus \Rr \Delta_d, \ | \
b_0(\Rcsi) \geq g(\csi) +1 - ad\} ; $$
it is an open cone in $\Rhd$. The main purpose of this
article is to prove the following
\begin{theo} \label{principal} Let $(X,c_X)$ be a  smooth real projective surface
and $(L,c_L)$ be a real Hermitian holomorphic line bundle on $X$ with
positive curvature. Then for every  sequence $d\in \Nn^* \mapsto a(d)\geq 1$ of rationals, 
there exist constants $C$, $D>0$, such that
$$ \forall d \in \Nn^*, \ \mu(\bm_d^{a(d)}) \leq Cd^6 e^{-D \frac{d}{a(d)}},$$
where $\mu(\bm_d^{a(d)}) $ denotes the Gaussian measure of
$\bm_d^{a(d)}$.
\end{theo}
The Gaussian measure $\mu$ on the Euclidian space
$(\Rhd, \langle \ , \  \rangle)$ is defined by the formula
$$ \forall A \subset \Rhd , \  \mu(A) =
\frac{1}{\sqrt \pi^{N_d} } \int_A e^{-|x|^2} dx,$$ where
$dx$ denotes the Lebesgue measure associated to $\langle \ , \ \rangle$.
This Gaussian measure is a probability measure on $\Rhd$
 invariant under its isometry group.

\begin{remark}
Likewise, the scalar product $\langle \ , \ \rangle$ induces a Fubini-Study form on the linear system $P(\Rr H^0(X,L^d))$. The volume of the projection $P(\bm^a_d)$ for the associated volume form 
just coincides with the measure $\mu(\bm^a_k)$
computed in Theorem \ref{principal}.
\end{remark}

In particular, when the sequence $a(d)$ is bounded, Theorem \ref{principal} implies
that the measure of the set $\bm^{a(d)}_d$ decreases exponentially 
with the degree $d$. This exponential rarefaction holds in particular for 
the set of maximal curves. 

\subsection{Real curves } 

When $X$ is one-dimensional, we get the following result: 
\begin{theo}\label{dimension1}
Let $X$ be a  closed  smooth real curve and $L$ be a real Hermitian holomorphic line bundle  on $X$ with
positive curvature. Then for every positive sequence $(\epsilon(d))_{d\in \Nn}$
of rationals numbers, there exist constants $C$, $D>0$ such that
$$ \forall d \in \Nn, \ \mu \{ \si \in \Rhd \setminus \Rr \Delta_d, \ | \
\# (\si^{-1}(0) \cap \Rr X) \geq \sqrt d \epsilon (d) \} \leq Cd^{3/2} e^{-D\epsilon^2(d)}, $$
where  $\mu$ denotes the Gaussian measure of the space $\Rhd$.
\end{theo}

\subsection{Roots of polynomials}  When $(X,c_X)$ is the projective
space of dimension $n\geq 1$, and $(L,h,c_L)$ is the degree one 
holomorphic line bundle equipped with its standard Fubini-Study metric, 
the vector space $\Rhd$ gets isomorphic to the space  $\Rr_d[X_1, \cdots, X_n]$
of polynomials with $n$ variables, real coefficients and  degree at most $d$. 
The scalar product induced on $\Rr_d[X_1, \cdots, X_n]$ by this isomorphism is the 
 one turning  the  basis  
$ \left(\sqrt {\binom {d+n}{j}} x_1^{j_1}\cdots x_n^{j_n}\right)_{0\leq j_1+ \cdots + j_n \leq d}$
into an orthonormal one,
where $  \binom {d+n}{j} = \frac{(d+n)!}{n!j_1!\cdots j_n! (d-j_1- \cdots - j_n)}.$
Thus, the induced measure $\mu$ on $\Rr_d[X_1, \cdots, X_n]$ is the Gaussian
measure associated to this basis. As a special case of our Theorems
\ref{principal} and \ref{dimension1}, we get the following
\begin{coro}
\begin{enumerate}

\item For every positive sequence $(\epsilon(d))_{d\in \Nn^*} $ of rational numbers,
 there exist positive constants $C$, $D$ such that the measure of the space of polynomials $P\in 
\Rr_d[X]$ which have at least $\epsilon (d) \sqrt d $ real
roots is bounded by $Cd^{3/2} \exp (-D\epsilon^2(d))$. 

\item For every sequence $d\in \Nn^* \mapsto a(d)\geq 1$ of rationals, there
exist positive constants $C,D$ such that the measure
of the space of polynomials $P\in  \Rr_d[X,Y]$ whose vanishing locus in $\Rr^2$ has at least $\frac{1}{2} d^2 - da(d)$ connected components
is bounded from above by $C d^6 \exp(-D\frac{d}{a(d)})$. 
\end{enumerate}
\end{coro}

\subsection{Strategy  of the proof of Theorems \ref{principal} and \ref{dimension1}} 
Every curve $C_\si \subset X$, $\si \in \Rhd \setminus \{0\}$, defines a current of integration 
which we renormalize by $d$, for  its mass 
not to depend on $d\in \Nn^*$. In order to prove Theorems \ref{principal}
and \ref{dimension1}, we first obtain  large deviation  estimates for
the random variable defined by this current. When $d$ grows to infinity, the expectation  of this variable converges
outside of $\Rr X$ to  the curvature form $\omega$ of $L$. 
These results thus go along the same lines as the one of 
Shiffman and Zelditch \cite{Shi-Zel}. 
They make use in particular of the asymptotic  isometry theorem of Tian \cite{Tian} (see also \cite{Catlin} and \cite{Zelditch}), as well as
smoothness results \cite{Kerzman} and behaviour close to the diagonal for the Bergman kernel  \cite{Shi-Zel2},\cite{Be-Be-Sz}. In order to deduce from these  results informations on the random
variable $b_0$, we use the theory of laminar currents
introduced by Bedford, Lyubich and Smillie \cite{Be-Ly-Sm}. Indeed, 
we show as a corollary of a theorem of de Th\'elin \cite{deThelin}
that every current in the closure of the ones arising from $\bm^a_d$ (in particular
every  limit current of a sequence of real maximal curves) is  weakly laminar 
outside of the real locus $\Rr X$, see Theorem \ref{laminaire}. As a consequence, these
currents remain in a compact set away from $\omega$. At this point,
our large deviation estimates provide the exponential decay.

\subsection{Description of the paper } In the first paragraph, we  bound the
Markov moments needed for our large deviation estimates. In the second
paragraph, we recall some elements of the theory of laminar currents in order to establish Theorem \ref{laminaire}, that is  laminarity
outside of the real locus of currents in the closure of the union of the sets
$\bm^a_d$.  We then get our estimates and their corollaries. We   prove   Theorems \ref{principal}
and \ref{dimension1} in the third paragraph, dealing separately with the cases of bounded and unbounded  sequences $a$. The last paragraph is devoted to some final discussion on the existence
of real maximal curves on general real projective surfaces as well 
as on the expectation of the current of integration on real divisors.

\textit{Acknowledgements.} We are grateful to  C\'edric Bernardin for 
 fruitful discussions on Markov moments and large deviations. This work was supported
 by the French Agence nationale de la recherche, ANR-08-BLAN-0291-02.

\section{Markov-like functions on real linear systems of divisors}

\subsection{Case of  projective spaces}\label{paragraph}
Let $\bo_{\Cc P^k} (1)$ be the degree one line bundle over
$\Cc P^k$ and $||. ||$ be its standard Hermitian
metric with curvature the Fubini-Study K\"ahler form $\omega_{FS}$. For every integer
$m\geq 1$, set

\beq \label{inequality}
M^m_{\Cc P^k}: \Cc P^k &\to& \Rr \\
z &\mapsto &\int_{\Rr H^0(\Cc P^k, \bo_{\Cc P^k} (1))}
                        \left|\log ||\si(z)||^2\right|^m d\mu(\si),
\eeq
where $d\mu$ denotes the Gaussian measure
of  the Euclidian  space $\Rr H^0(\Cc P^k, \bo_{\Cc P^k} (1))$.

\begin{prop}\label{tau}
For every $m\geq 1$, the function $M^m_{\Cc P^k}$ satisfies
$$ \forall z \in \Cc P^k\setminus \Rr P^k, \  M^m_{\Cc P^k}(z) \leq \frac{4 m! (k+1)}{1-||\tau(z)||},$$
where $\tau$ denotes the section of $\bo_{\Cc P^k}(2)$ defined
by $\tau(z_0, \cdots, z_k) = z_0^2 + \cdots  + z_k^2$.
\end{prop}
\begin{remark}\label{invariance}  The holomorphic section $\tau$
in Proposition \ref{tau} is invariant under the action of the group
$PO_{k+1}(\Rr)$ of real isometries of $\Cc P^k$.
A slice   of $\Cc P^k$ for this action
is given by the interval $ I = \{z_r = [1: ir: 0 : \cdots : 0],   0\leq r\leq 1\},$
where the end $ r= 0$ (resp. $r= 1$) corresponds
to the orbit $\Rr P^k$ (resp. $\tau^{-1}(0)$) of this action.
\end{remark}

\bpr [ of Proposition \ref{tau}] Both members of the inequality are invariants
under the action of $PO_{k+1}(\Rr)$, so that it
is enough to prove it for $ z$ in the fundamental domain $I$.
Let $\si_0, \cdots, \si_k$ be
the  orthonormal basis of $\Rr H^0(\Cc P^k, \bo_{\Cc P^k} (1))$
given by $\si_i ([z_0: \cdots : z_k]) = \sqrt{k+1} z_i.$
This basis induces the isometry
$$ a = (a_0, \cdots , a_k) \in \Rr^{k+1} \mapsto
\si_a = a_0\si_0+ \cdots + a_k\si_k \in \Rr H^0(\Cc P^k, \bo_{\Cc P^k} (1)).$$
By definition, for every $a \in \Rr^{k+1}$ and $z\in \Cc P^k$,
$$ ||\si_a(z) ||^2 = (k+1) \frac{|a_0 z_0+ \cdots + a_k z_k|^2}{|z|^2}.$$
We deduce that for every $0< r\leq 1$ and $m\geq 1$,
\beq
M^m_{\Cc P^k}(z_r) &=& \int_{\Rr^{k+1}}  \left|\log ||\si_a (z_r)||^2\right|^m d\mu(a)\\
& =&  \int_{\Rr^{2}} \left| \log  \left((k+1) \frac{|a_0 + ira_1|^2}{1+r^2} \right)\right|^m \frac{e^{-|a|^2}}{\pi} da_0da_1\\
& = & \int^{\infty}_0 \int_0^{2\pi} \left| \log  ((k+1)\rho^2) + \log  \frac{|\cos \theta + ir \sin \theta|^2}{1+r^2}  \right|^m
\frac{e^{-\rho^2}}{\pi} \rho d\rho d\theta\\
& \leq & 2 \int_0^\infty \left| | \log ((k+1) \rho^2) |+\log (\frac{1+r^2}{r^2}) \right|^m e^{-\rho^2}\rho d\rho,
\eeq
since for every $0\leq r\leq 1$ and every $  \theta \in [0,2\pi],$
$ r^2 \leq \cos^2 \theta + r^2 \sin^2\theta =
| \cos \theta + i r \sin \theta|^2 \leq 1. $
Hence,
$$
M^m_{\Cc P^k}(z_r) \leq 2 \int_0^1 (-\log(\frac{\rho}{\alpha})^2)^m e^{-\rho^2} \rho d\rho
+
2 \int_1^\infty (\log(\alpha \rho)^2)^m e^{-\rho^2} \rho d\rho,$$
where $$\displaystyle \alpha^2 = (k+1) \frac{(1+r^2)}{r^2} = 2\frac{(k+1)}{1 - ||\tau (z_r)||}. $$
We  now compute  these two integrals.
\beq  
         2                  \int_0^1 (-\log(\frac{\rho}{\alpha})^2)^m e^{-\rho^2} \rho d\rho
&= &2\alpha^2      \int_0^{1/\alpha}  (-\log \rho^2)^m e^{-\alpha ^2 \rho^2} \rho d\rho \\
&\leq & 2\alpha^2 \int_0^1  (-\log \rho^2)^m  \rho d\rho \\
&=& \alpha^2 \int_0^\infty t^m e^{-t} dt \ \ \text{with  } \ t= -\log \rho^2   \\
&= &2 \frac{(k+1) m! }{1 - ||\tau(z_r)||}.
\eeq
As for the second integral, we deduce from the estimate  $\rho e^{-\rho ^2} \leq e^{-1/\rho^2}/\rho^3 $ valid for every  $\rho\geq 1$, that
$$ 2 \int_1^\infty (\log(\alpha \rho)^2)^m e^{-\rho^2} \rho d\rho
\leq 2 \int_1^\infty (\log(\alpha\rho)^2)^m \frac{e^{-1/\rho^2}}{\rho^3} d\rho.$$
With $t = 1/\rho$, the right hand side  becomes  
$\disp 2\int_0^1  (-\log(\frac{t}{\alpha})^2)^m e^{-t^2} t dt \leq  2 \frac{(k+1) m! }{1 - ||\tau(z)||}.$
\epr

\subsection{Asymptotic results in the general case}\label{asymptotic}
Let $X$ be a closed real  K\"ahler manifold of dimension $n$ and
$ L$ be a real Hermitian  line bundle  over $X$ with positive curvature $\omega$.
We denote by $d_L$ the smallest integer such that $L^{d}$ is very ample
 for every $d\geq d_L$ and by $ \Phi_d : X \to P (H^0(X,L^d)^*)$
the associated  embedding, where $x\in X$ is mapped to
the set of linear forms that vanish on the hyperplane $\{\si \in H^0(X,L^d) \ | \ \si (x) = 0\}$.
The $L^2$-Hermitian product of $\hd$ induces a Fubini-Study metric
on the complex projective space $P(\hd^*)$ together with a Hermitian metric $|| . ||$ on
the bundle $\bo_{P(\hd^*)} (1)$. We denote by  $h_{\Phi_d}$ the pullback of  $||  . ||$ 
on $L^d$ under the canonical isomorphism
$ L^d \to \Phi_d ^* \bo_{P(\hd^*)} (1). $
Recall that the latter is induced by the isomorphism
$$ (x,\alpha)\in L^d \mapsto (\Phi_d(x), \alpha_x) \in  \bo_{P(\hd^*)} (1), $$
where $\alpha_x : \si \in H^0(X, L^d) \mapsto \langle \sigma, \alpha \rangle_x \in \Cc$.
In particular, the curvature form of $h_{\Phi_d}$ is $\Phi_d ^* \omega_{FS}$, where $\omega_{FS}$
denotes the Fubini-Study form of $P (\hd^*)$.
The quotient $h^d/h_{\Phi_d}$ of these metrics of $L^d$ is given by
the function
$ x\in X \mapsto \sum_{i=1}^{N_d} h^d(\si_i(x), \si_i(x)),$ where
$(\si_1, \cdots , \si_{N_d})$ stands for any orthonormal basis of $\hd$.
Let $|| . ||_{\Phi_d}$ be the norm induced by $h_{\Phi_d}$.
For every  $m \in \Nn^*$, we set
\beq 
M^m_{(X,L^d)} : X &\to &\Rr \\
x & \mapsto & \int_{\Rr \hd}
                        \left|\log ||\si(z)||_{\Phi_d}^2\right|^m d\mu(\si),
\eeq
so that $M^m_{(X,L^d)} = M^m_{P(\hd^*)} \circ \Phi_d$.

\begin{prop}\label{sup}
Let $X$ be a closed real  K\"ahler manifold of dimension $n$ and
$ L$ be a positive real Hermitian  line bundle  over $X$.
For every sequence $(K_d)_{d\in \Nn^*}$ of compact subsets 
of $X\setminus \Rr X$ such that the sequence $(d\text{ dist} (K_d, \Rr X)^2)_{d\in \Nn^*}$
grows to infinity, the sequence $ \left| \left| \ || \tau|| \circ \Phi_d  \right| \right| _{C^0(K_d)}$
converges to zero as $d$ grows  to infinity, where $\tau \in \Rr H^0 (\bo_{ P(\hd^*)}(2))$ is the section defined in Proposition \ref{tau}. If $K$ is a fixed compact, the same holds for any norm $C^q(K)$, $q\in \Nn$.
\end{prop}
The $L^2$-Hermitian   product on $\hd$ induces a scalar product
on $ \Rhd$ and its dual $\Rhd^*$. Let $\langle \ , \ \rangle_\Cc$ be the extension
of this scalar product to a  complex bilinear product on $\hd^*$. 
The section $\tau$ of $\bo(2)_{P(\hd)^*}$ that appears in 
Propositions \ref{tau} and  \ref{sup}
is the one induced by $ \si^* \in \hd^* \mapsto \langle \si^*, \si^* \rangle_\Cc \in \Cc.$

\bpr [ of Proposition \ref{sup}]
Let $D^*\in L^*$ be the open unit disc bundle for the metric $h$ and
$dx dv$ the product measure on $D^*$, where $ dx = \omega^n/\int_X \omega^n$
is the measure on the base $ X$ of $D^*$ and $dv$ is the Lebesgue measure
on the fibres. Let $L^2(D^*, \Cc)$ be the  space of complex functions of class $L^2$ on $D^*$
for the measure $dx dv$ and $\bh^2 \subset L^2(D^*,\Cc)$ be the closed subspace
of $ L^2$ holomorphic functions on the interior of $D^*$.
Every function $f \in \bh^2$ has a unique expansion
$$ f : (x,v) \in D^* \mapsto \sum_{d=0}^{\infty} a_d(x)v^d \in \Cc,$$
where for every  $d\geq 0$, $a_d \in H^0(X,L^d)$. The series
$ \disp \sum_{d=0}^{\infty} a_d(x)v^d$ converges uniformly on every compact subset of $D^*$
as well as  in $L^2$-norm.
Let $B$ be the Bergman kernel of $D^*$, defined
by the relation:
$$ \forall (y,w) \in D^*, \forall f\in \bh^2, \ f(y,w) = \int_{D^*} f(x,v) \overline{B((x,v),(y,w))} dxdv,$$
where this function $B : D^* \times D^* \to \Cc$
is holomorphic in the first variable and antiholomorphic in the second one.
Kerzman \cite{Kerzman} proved that the Bergman kernel can be extended smoothly  up to
the boundary outside of the diagonal of $\overline{D^*}\times  \overline{D^*} $.
Now, denote by  $ D^*_{\RX} = \{(x,v) \in D^* \ | \ x\in \RX \}. $
The function 
$$ b : (x,v) \in D^*_{\RX} \mapsto B((x,v), c_{L^*}(x,v)) \in \Cc $$
is holomorphic on $D^*_{\RX}$ and can be extended smoothly on $\overline{D^*_{\RX}}$. 
For every $d \geq 0$ let $(\si_{i,d})_{1\leq i\leq N_d}$ be an orthonormal basis
of $\Rhd$ and $$e_{i,d} = \sqrt{\frac{d+1}{\pi}} \widehat {\si_{i,d}} \in \bh^2,$$
where $\widehat {\si_{i,d}} : (x,v) \in D^* \mapsto \si_{i,d}(x)v^d\in \Cc$.
The family  $(e_{i,d})_{i,d}$  forms a Hilbertian basis
of $\mathcal H^2$. 
Then,
\beq
\forall (x,v) \in D^*_{\RX}, \ b(x,v) &=& B( (x,v), c_{L^*}(x,v) )\\
& = & \sum_{d=0}^{\infty} \sum_{i=1}^{N_d} e_{i,d}(x,v) \overline{e_{i,d}(c_{L^*}(x,v))}\\
& = &  \sum_{d=0}^{\infty} \frac{d+1}{\pi} \sum_{i=1}^{N_d} \si_{i,d}(x)v^d \overline{\si_{i,d}(c(x)) c_{L^*}(v)^d}      \\
& = &  \sum_{d=0}^{\infty} \frac{d+1}{\pi} \sum_{i=1}^{N_d} \si^2_{i,d}(x)v^{2d}
\eeq
since $\si_{i,d}$ is real, that is satisfies $ \si_{i,d}\circ c = c_{L^d}\circ \si_{i,d}$.
Note that $ \tau\circ \Phi_d = \sum_{i=1}^{N_d} \si^2_{i,d}(x) $ and
from  Tian's asymptotic isometry theorem \cite {Tian} (see also \cite{Catlin} and \cite{Zelditch} )
$ \disp ||. ||_{\Phi_d} \leq \frac{C}{d^n}|| . ||$. For a fixed compact $K$,
the result thus just follows from Cauchy formula applied to $b$. 
In general, we may substitute $L$ with $L^d$, $K$ with $K_d$ and deduce
the result from Proposition 2.1 of \cite{Shi-Zel2} (see also \cite{Be-Be-Sz}) since the sequence $(d \sup_{x\in K_d} \text{ dist}(x,c_X(x))^2)_{d\in \Nn^*}$ grows to infinity  as
$d$ grows to infinity.
\epr

Proposition \ref{sup} and Proposition \ref{tau} imply the following  
\begin{coro}\label{corollaire} Let $X$ be a real  K\"ahler manifold and
$ L$ be a real positive Hermitian  line bundle  over $X$.
For every sequence $(K_d)_{d\in \Nn}$ of compact subsets of $\RX$ such that the 
sequence $(d \text{ dist}(K_d, \Rr X)^2)_{d\in \Nn^*}$ grows to infinity, there exists
a positive constant $c_K$ such that as soon as $L^d$ is very ample, 
  $$ \forall m\in \Nn^*,  \sup_{K_d} M^m_{(X,L^d)} \leq c_K m! N_d,$$
where $N_d = \dim \hd$.
\end{coro}

\begin{remark}\label{exp} 
Actually, in Corollary \ref{corollaire},  $\lim_{d\to \infty} 1/N_d \sup_{K_d} M^m_{(X,L^d)} \leq 4m!$. 
Also, Proposition \ref{sup} and even the exponential decay of the quantity $\underset{K}{\sup} || \tau|| \circ \Phi_d$
are easy to establish in some cases, including the following ones. 
\end{remark}

\textit{Projective spaces}. When $X = \Cc P^n$, $\Phi_d : \Cc P^n \to \Cc P^{N_d-1}$ is equivariant with respect to the
groups of real isometries $PO_{n+1}(\Rr)$ and $ PO_{N_d}(\Rr)$. Since $\tau$
is invariant  under these actions, $\tau\circ \Phi_d$ has to be a multiple of the section $\tau^d$. Now $|| \tau||_{\Rr P^n} \equiv 1$,
so that $ \tau \circ \Phi_d = \tau^d$ and 
$\underset{K}{\sup}  || \tau \circ \Phi_d|| = \underset{K}{\sup}  || \tau ||^d .$
This was observed by Macdonald \cite{Macdonald} in the case  $X = \Cc^n$. \\

 \textit{Ellipsoid quadrics.}
Assume now that  $X= \{[x_0: \cdots : x_{n+1}]\in \Cc P^{n+1} \  | \ x_0^2 = x_1^2 +\cdots + x_{n+1}^2          \} $
is the ellipsoid quadric and $L$ the restriction of $\bo_ {\Cc P ^{n+1}}(1)$ to $X$. Then,
$\tau \circ \Phi_d$ is a multiple of the hyperplane section $x_0^{2d}$,
since it is invariant under the group of isometries $O_{n+1}(\Rr)$ acting
on the coordinates $(x_1, \cdots, x_{n+1}). $ At a real point $x = (x_0, \cdots , x_{n+1})$,
 $|| \tau ||\circ \Phi_d = 1 $ and
$$ ||x_0^2|| =\frac{|x_0^2|}{|x_0|^2 + \cdots + |x_{n+1}|^2} = \frac{1}{2}
\left( \frac{x_0^2}{x_0^2 + \cdots + x_{n+1}^2}+
\frac{x_1^2 + \cdots + x_{n+1}^2}{x_0^2 + \cdots + x_{n+1}^2}\right) = \frac{1}{2},$$
hence $\tau\circ \Phi_d= 2^d x_0^{2d}$ and
$\underset{K}{ \sup} ||\si|| \circ \Phi_d = \underset{K}{ \sup} (2||x_0^2||)^d.$

\textit{The hyperboloid surface.}
Finally, if $X  = (\Cc P^1_1\times \Cc P^1_2; \conj \times \conj)$ is the hyperboloid
quadric  surface
and $L = \bo(a)_{ \Cc P^1_1} \otimes  \bo(b)_{ \Cc P^1_2} $ with $a,b>0$,
then $\tau\circ \Phi_d$ is $PO_2(\Rr)\times PO_2(\Rr)$-invariant,
hence a multiple of $(\tau_1^a\otimes \tau_2^b)^d$ where
$\tau_i =\tau _{ \Cc P^1_i}\in \bo(2)$ for $i=1,2$. Computed at a real point,
this multiple is one, so that $\tau\circ \Phi_d = (\tau_1^a\otimes \tau_2^b)^d$ and
$\underset{K}{ \sup} ||\si|| \circ \Phi_d = (\underset{K}{ \sup} ||\tau_1^a\otimes \tau_2^b||)^d.$

\section{Weakly laminar currents and large deviation estimates} 
\subsection{Weakly laminar currents} \label{paragraph2}
Let $(X,\omega)$ be a smooth K\"ahler manifold and $\bt^{(1,1)}_{L^2} $ be its space of closed positive currents of type $(1,1)$
and  mass $ L^2 = \int_X \omega \wedge \omega$. 
Recall that by definition such a current is a continuous linear form 
on the space of smooth two-forms that
vanishes on  forms of type $(2,0)$ and $(0,2)$ as well as on the exact
forms. Moreover, the mass $\langle T,\omega \rangle $ equals
 $L^2$ and $T$ is positive
once evaluated on a positive $(1,1)$-form. In particular, 
$T$ is of measure type, that is continuous for the sup norm on 
two-forms. The space $\bt^{(1,1)}_{L^2} $, equipped with the weak topology,
is a compact and convex space. 
For every $d\in \Nn^*$ and every $\si\in \hd\setminus\{0\}$, denote by 
$Z_{\si} \in \bt^{(1,1)}_{L^2}$ the current of integration 
$$ Z_\si : \phi \in \Omega^2(X) \mapsto \frac{1}{d}\int_{C_\si} \phi \in \Rr,$$
where $C_\si = \si^{-1}(0)$. 
The following definition of weakly laminar currents was introduced in \cite{Be-Ly-Sm}.
\begin{defi}\label{wl} A current $T \in \bt^{(1,1)}_{L^2}$ is called weakly laminar in the open set $U\subset X$ 
iff there exist a family of embedded discs $(D_a)_{a\in \ba}$ in $U$ together with a measure $da$ on 
$\ba$, such that for every $a$ and $a'$ in 
$\ba$, $D_a\cap D_{a'}$ is open in $D_a$ and $D_{a'}$, and such
that for every smooth two-form $\phi$ with support in $U$, 
$$ \langle T,\phi \rangle = \int_{a\in \ba}\left(\int_{D_a} \phi \right) da.$$
\end{defi}

For every open subset $U$ of $X$, denote by  $Lam(U)\subset \bt^{(1,1)}_{L^2}$
the subspace of  closed  positive 
currents of mass $L^2$ which are weakly laminar on $U$. For
all $a\in \Qq^*_+$, denote by $\bz^a$ the closure of the union $\cup_{d\in \Nn^*} \bz_d^a$
in $  \bt^{(1,1)}_{L^2}$, where $ \bz_d^a $ denotes the image of the set  
$$ \bm^a_d = \{\si \in \Rhd \setminus \Rr \Delta_d \ |Ê\ b_0(\Rr C_\si)\geq g(C_\si) + 1 - ad\}$$
under the map $\si \mapsto Z_\si \in \bt^{(1,1)}_{L^2}$. 
\begin{theo}\label{laminaire}
Let $X$ be a closed real  projective surface and $L$ be a positive real Hermitian  line bundle on $X$. Then, for every $a\in \Qq^*_+$, the inclusion $\bz^a \subset Lam(X\setminus \Rr X) $ holds.
\end{theo}
In particular, every limit of a sequence of real maximal curves is weakly laminar
outside of the real locus of the manifold.

\bpr
This result is actually a direct consequence of    Theorem 1 of \cite{deThelin}.
Indeed, let $T\in \bz^a$ and $(Z_{\si_d})_{d\in \Nn^*}$  be a sequence of 
currents of integration which converges to $T$ (we can
indeed assume that $ T\notin \bigcup_{d\in \Nn^*}  \bz_d^a$). 
For every $d\in \Nn^*$, the genus of the complement $C_{\si_d} \setminus  \Rr C_{\si_d}$
satisfies $g(C_{\si_d} \setminus  \Rr C_{\si_d})\leq ad$, while 
its area equals $d\int_X \omega^2$. Let $B$ be a ball
with compact closure in $ X\setminus \Rr X$ and $A(C_{\si_d}\cap B)$ be the area of the restriction of 
$C_{\si_d}$ to $B$ . Without loss of generality, we can assume that $\frac{1}{d} A(C_{\si_d}\cap B)$
converges to $m_B \in [0, \int_X \omega^2]$. If $m_B= 0$, the restriction of $T$ to $B$ vanishes.
Otherwise, $g(C_{\si_d} \cap B) = O( A(C_{\si_d}\cap B)),$
where  the area $A(C_{\si_d}\cap B)$ can be computed for the flat metric on the ball. 
From Theorem 1 of \cite{deThelin} we know that $1/m_B  T_{|B}$
is weakly laminar. Hence the result.
\epr
\begin{lemme}\label{lemme} Let $\omega$ be  a K\"ahler form on a complex surface $X$. Then $\omega$ is 
 nowhere weakly laminar.
\end{lemme}

\bpr
Assume that there exists an open subset $U$ of $X$ and a measured family $(D_a)_{a\in \ba}$
of embedded discs in $U$ given by Definition \ref{wl}  such that for every two-form $\phi$
with compact support in $U$, 
$$ \int_U \omega\wedge \phi = \int_\ba \left(\int_{D_a} \phi \right) da. $$
For every two-form $\psi$
defined and continuous on $ \bigcup_{a\in \ba} D_a$, we denote by $T_\psi$ the current 
$$ \phi \in \Omega^2_c(U) \mapsto T_\psi(\phi) = 
\int_{a\in \ba} \left( \int_{D_a} (\frac{\psi\wedge \phi}{\omega^2})\omega\right) da.$$
Then $T_\omega = \omega$, since for every $ \phi \in \Omega^2_c(U)$, 
$$ T_\omega(\phi) =
\int_{ \ba} \left( \int_{D_a} (\frac{\omega\wedge \phi}{\omega^2})\omega\right) da
= \int_U (\frac{\omega\wedge \phi}{\omega^2})\omega^2 = \int_U \omega \wedge \phi.
$$
Let $\psi$ be the (1,1)-form  defined along $\disp \cup_{a\in \ba} D_a$ 
in such a way that for every $a\in \ba$ and  $x\in D_a$, $T_xD_a$ lies in the kernel of
$\psi_x$ and $\psi_x\wedge \omega = \omega^2$. Then  
$T_\psi = \omega$, since
$$ \forall \phi \in \Omega^2_c(U), \ T_\psi(\phi) = 
\int_{ \ba} \left( \int_{D_a} (\frac{\psi\wedge \phi}{\omega^2})\omega\right) da = 
 \int_{ \ba} \left( \int_{D_a}  \phi\right) da.$$
But $\psi$ is integrable for both currents $T_\psi$ and $T_\omega$ 
and $T_\psi(\psi) = 0$ while $T_\omega(\psi) =  \omega^2$. This
contradicts the equality $T_\psi =T_\omega$. 
\epr

\subsection{Large deviation estimates}
\begin{prop}\label{prop}
Let $X$ be a real projective manifold of dimension $n$ and $L$
an ample real Hermitian line bundle over $X$. For every smooth 
$(2n-2)-$form $\phi$ with compact support  in  $X \setminus \Rr X$,
every  $ d\geq d_L $ and every $ \epsilon >0$, the measure of  the set 
$$  \left\{ \si \in \Rhd\setminus \Rr \Delta_d \ | \  \frac{1}{d} \left|\int_X \log ||\si(x) ||^2_{\Phi_d} \partial \del \phi \right| \geq 
\epsilon \right\}  $$
is bounded from above by the quantity 
$$  \frac{2 c_{K_\phi} N_d}{Vol(K_\phi)} 
\exp \left(\frac{-\epsilon d}{2||\partial \del \phi||_{L^\infty} \text{Vol }(K_ \phi)}\right),$$
where $K_\phi$ can denote both the support of $\phi$ or
the support of $\partial \del \phi$, while $c_{K_ \phi}$ and $d_L$ are given in \textsection \ref{asymptotic}.
\end{prop}
Recall that $X$ is equipped with the volume form $dx = \omega^n/\int_X \omega^n$. 
The norm $ ||\partial \del \phi||_{L^\infty}$ and the volume $Vol(K_\phi)$ are computed with respect to this volume form. Recall also that $N_d$ denotes the dimension of $\Rhd$ and finally that from the Poincar\'e-Lelong formula, the current 
$ 1/(2i\pi d) \partial\del \log ||\si(x) ||^2_{\Phi_d} $ coincides with 
$ \frac{1}{d}\Phi_d^* \omega_{FS} - Z_\si.$

\bpr We use Markov's trick. For every $ \lambda >0$,
$$ \left| \int_X \log ||Ê\si||^2_{\Phi_d} \partial \del \phi \right| \geq 
 \epsilon d \Longleftrightarrow
  \exp \left( \lambda \left|\int_X  \log ||\si||^2_{\Phi_d} \partial \del \phi \right|\right)
 \geq e^{\lambda  \epsilon d}, $$
where
\beq
 \exp \left(\lambda \left|\int_X  \log ||\si||^2_{\Phi_d} \partial \del \phi \right| \right) &=& 
 \sum_{n=0}^{\infty} \frac{\lambda^m}{m!}  \left| \int_X  \log ||\si||_{\Phi_d}^2 \partial \del \phi \right|^m. 
 \eeq
 From H\"older's inequality we get
\beq
 \left| \int_X  \log ||\si||_{\Phi_d}^2 \partial \del \phi \right|^m &  \leq &
                               \int_X \left| \log ||\si||_{\Phi_d}^2 \right|^m dx     
               \left(          \int_X  \left| \partial \del \phi \right|^\frac{m}{m-1}  dx  \right)^{m-1} \\
               & \leq & \frac{1}{Vol(K_\phi)} \left(|| \partial \del \phi ||_{L^\infty}
               Vol (K_\phi)\right)^m   \int_X \left| \log ||\si||_{\Phi_d}^2 \right|^m dx.
 \eeq
As a consequence, for every $d\geq d_L$, 
the measure $\mu_\epsilon^d(\phi)$ of our set satisfies
\beq 
e^{\lambda \epsilon d}Ê\mu_\epsilon^d(\phi) &\leq &
        \int_{\Rhd} \exp \left(        \lambda \left|  \int_X  \log ||\si||^2_{\Phi_d} \partial \del \phi \ \right| d\mu(\si)\right)\\
& \leq &  \frac{1}{ Vol(K_\phi ) } \sum_{n=0}^{\infty} \frac{\lambda^m}{m!} || \partial \del \phi ||_{L^\infty}^m
               Vol (K _ \phi)^m   \int_X M^m_{(X,L^d)} dx, 
            \eeq
            where $M^m_{(X,L^d)}$ is defined in \textsection 1.2. 
            Thanks to Corollary \ref{corollaire}, the latter right hand side is bounded from above by               
$\frac{c_{K_\phi }N_d}{Vol (K_\phi)}  \sum_{n=0}^{\infty} (\lambda  ||\partial \del \phi ||_{L^\infty} Vol (K_\phi))^m 
$, that is  $$ \frac{c_{K_\phi} N_d}{Vol (K_\phi)(1- \lambda || \partial \del \phi ||_{L^\infty} Vol (K_\phi))}.$$
The result follows by choosing $ \lambda = (2|| \partial \del \phi ||_{L^\infty} Vol (K_\phi))^{-1}.$
\epr
\begin{coro}\label{corollaire2}
Under the hypotheses of Proposition \ref{prop}, let $\stackrel{\circ}{K}$ be
a relatively compact open subset of $X \setminus \Rr X$. Then, there exist
constants $C_K$, $D_K$, $\lambda_K>0$ and, for every $d\geq d_L$, 
a subset $\ba^d_K \subset \Rhd$ of measure bounded from above by $ C_K e^{-D_Kd }$
such that for every $\si \in \Rhd \setminus \ba^d_K$, 
the volume $A(C_\si\cap K)$ of $C_\si\cap K$ satisfies $A(C_\si\cap K) \geq \lambda_K d$. 
\end{coro}

\bpr Let $\stackrel{\circ}{K_1}$ be a relatively compact open subset of $\stackrel{\circ}{K}$ and 
$\chi : X  \to [0,1]$ be a smooth cutoff function 
with support in $K$ such that $\chi_{|K_1} \equiv 1$. 
Applied to $\phi = \chi \omega^{n-1}$
and $\epsilon = \pi \int_{K_1} \omega^n$, Proposition \ref{prop} 
provides us with constants $C_K$ and $D_K$ such that the set 
$$ \ba^d_K = \{0 \}\cup \left \{ \si \in \Rhd\setminus \{0\}   \ | \
\frac{1}{d} \left| \int_X  \log ||\si ||^2_{\Phi_d} \partial \del (\chi\omega^{n-1}) \right |    
\geq \pi \int_{K_1} \omega^n   \right \}$$
is of measure bounded from above by  $C_Ke^{-D_K d}$, 
since $N_d$ grows polynomially with $d$.
From Poincar\'e-Lelong formula follows that  for every $\si \in \Rhd \setminus \ba^d_K$, 
$$ \left|  \frac{1}{d}\int_{C_\si} \chi \omega^{n-1} - \int_X \frac{1}{d} \Phi_d^* \omega_{FS} \wedge \chi \omega^{n-1}\right| \leq \frac{1}{2}\int_{K_1} \omega^n.  $$
Now, from Tian's asymptotic isometry theorem \cite{Tian} (see also  \cite{Catlin} and \cite{Zelditch}),  
$\frac{1}{d}\Phi_d^*\omega_{FS}$ converges to $\omega$ as $d$
grows to $\infty$, so that for $d$ large enough,
$  \int_X \frac{1}{d} \Phi_d^* \omega_{FS} \wedge \chi \omega^{n-1} \geq \int_{K_1} \omega^n$. 
As a consequence, $$\frac{(n-1)!}{d} A(C_\si \cap K) \geq \frac{1}{d}\int_{C_\si} \chi \omega^{n-1}
\geq \int_X\frac{1}{d} \Phi_d^* \omega_{FS} \wedge \chi \omega^{n-1} - \frac{1}{2} \int_{K_1} \omega^n.$$ 
The left hand side being bounded from below by a constant $(n-1)! \lambda_K$, 
the result follows.
\epr
For every $n\in \Nn^*$ and $\rho>0$, denote by $B^{2n}(\rho)\subset \Cc^n$ the closed ball of radius 
$\rho$  and volume  $ \pi^n \rho^{2n}/n!$. The standard K\"ahler form of $\Cc^n$ 
is denoted by $\omega_0$.
\begin{defi}\label{ball}
Let $(X,\omega)$ be a K\"ahler manifold of dimension $n$. By abuse, we define
 a \textit{ball of radius $\rho$} of $X$ to be the image of a holomorphic
 embedding $\psi_\rho:  B^{2n}(\rho)\to X$ whose differential at the origin
 is  isometric and which everywhere 
 satisfies the inequalities $ 1/2 \omega_0 \leq \psi_\rho^* \omega \leq 2\omega_0$.
\end{defi}
\begin{coro}\label{corollaire3}
Under the hypotheses of Proposition \ref{prop}, let $\stackrel{\circ}{K}$ be 
a relatively compact  open subset of $X \setminus \Rr X$. Then, there exist
constants $D_K$, $\lambda^1_K$, $\lambda^2_K>0$
such that for every ball $B$ of radius $\rho>0$ included in $K$
 and  every $d\geq d_L$, there exists 
a set $\ba^d_B \subset \Rhd$ of measure 
$$\mu(\ba^d_B) \leq \frac {2c_K N_d \int_X \omega^n }{ \rho^{2n} } e^{-D_K d\rho^2 }$$
such that for every $\si \in \Rhd \setminus \ba^d_B$, 
the volume of $C_\si\cap B$ satisfies $\lambda^1_K d\rho^{2n}\leq A(C_\si\cap B) \leq \lambda_K^2 d\rho^{2n}.$
\end{coro}
Recall that the constant $c_K$ is given by Corollary \ref{corollaire}, while $d_L$
is defined in \textsection \ref{asymptotic}.

\bpr
Let $\chi : \Cc^n \to [0,1]$ be a smooth cutoff function 
with support in the unit ball and such that $\chi^{-1}(1)$ 
 contains the ball
of radius $\sqrt{2/\pi}$.
For any $\rho>0$, let $ \chi_\rho : x\in \Cc^n \mapsto \chi(x/\rho) $
be the associated cutoff function with support in the
ball of radius $\rho$. Let $\psi : B \to B^{2n}(\rho)\subset \Cc^n$ 
be a biholomorphism given by Definition \ref{ball},  
and $B_1 = (\chi_\rho \circ \psi)^{-1} (1)$. 
Let $\phi = (\chi_\rho \circ \psi)\omega^{n-1}$, $K_\phi = supp(\phi)$, and for every $d\geq d_L$,
$$ \ba^d_B = \{0 \}\cup \left \{ \si \in \Rhd\setminus \{0\}   \ | \
\frac{1}{d} \left|\int_X \log ||\si ||^2_{\Phi_d} \partial \del \phi \right|    
\geq \pi \int_{B_1} \omega^n   \right\}.$$
From Proposition \ref{prop} follows that
$$ \mu(\ba^d_B) \leq \frac{2c_{K_\phi} N_d}{ Vol (K_\phi)}
    \exp \left( \frac{-\pi d\int_{B_1} \omega^n}{2 ||\partial \del \phi ||_{L^\infty}     Vol (K_\phi)} \right),$$
    where $$(\int_X \omega^n )Vol(K_\phi) = \int_{K_\phi} \omega^n \geq \frac{1}{2^n} \int_{\chi_\rho^{-1}(1)} \omega_0^n
    \geq  \rho^{2n}.$$ 
    The metric on $B$ is bounded from above and below by the flat metric, see Definition \ref{ball}. 
    The quotient $\int_{B_1} \omega^n/Vol (K_\phi)$
    is thus bounded from below by a positive constant, since 
    $\int_{\chi_\rho^{-1}(1)}   \omega_0^n/ \int_{supp (\chi_\rho)} \omega_0^n$
    does not depend on $\rho$. Likewise, $|| \partial \del \phi||_{L^\infty}$ 
    is bounded from above by a multiple of 
    $ \sup_{B^{2n}(\rho)} \left| \partial \del \chi_\rho \wedge \omega^{n-1}_0/\omega^n_0 \right|$. The latter being of the order of $1/\rho^2$, 
     we  deduce the existence of a positive constant $D_K$ such
     that $$\mu(\ba^d_B) \leq 2(\int_X \omega^n) \frac{c_K N_d}{\rho^{2n}} exp (-D_K \rho^2 d).$$
    But for every $\si \in \Rhd\setminus \ba^d_B$, we have
    $$ \left | \frac{1}{d}\int_{C_\si} (\chi_\rho \circ \psi) \omega^{n-1} - \int_X\frac{1}{d}  \Phi_d^* \omega_{FS} \wedge
    (\chi_\rho \circ \psi) \omega^{n-1} \right | \leq \frac{1}{2}\int_{B_1} \omega^n.$$ 
    The term $\frac{1}{d} \int_X \Phi_d^* \omega_{FS} \wedge
    (\chi_\rho \circ \psi) \omega^{n-1}$ is greater than $$
     \int_{B_1} \omega^n + \int_{B\setminus B_1}  (\chi_\rho \circ \psi) \omega^{n-1}
     - || (\frac{1}{d} \Phi_d^* \omega_{FS}-\omega) \wedge \omega^{n-1}||_{L^\infty} Vol (B).$$
     From  Tian's asymptotic isometry theorem \cite{Tian}, $\frac{1}{d} \Phi_d^* \omega_{FS} $ converges to 
    $\omega$ as $d$ grows to infinity. Together with
    Definition \ref{ball}, this implies that
     for $d$ large enough, $$\frac{1}{d} \int_X \Phi_d^* \omega_{FS} \wedge
    (\chi_\rho \circ \psi) \omega^{n-1} \geq \int_{B_1} \omega^n$$
    and $$ \frac{(n-1) ! }{d\rho^{2n}}A(C_\si \cap B) \geq \frac{1}{d\rho^{2n}} \int_{C_\si} (\chi_\rho\circ \psi) 
    \omega^{n-1} \geq \frac{1}{2\rho^{2n}}\int_{B_1} \omega^n. $$
    The right hand side being bounded from below by a positive constant, we deduce the lower bound for $A(C_\si \cap B)$. Likewise, we deduce that 
    $  (n-1) ! /(d\rho^{2n}) A(C_\si \cap B) \leq 3/(2\rho^{2n})\int_{B} \omega^n$. 
    The right hand side being bounded from above by a positive constant, we deduce the 
    upper bound for $A(C_\si \cap B)$ replacing $B_1$ by $B$ in the proof.     
    \epr

 \begin{lemme}\label{lemme2}
 For every compact subset $K$ of a $n$-dimensional K\"ahler manifold, there 
 exist constants $r_K$ and $n_K>0$ such that for every $\rho>0$ small enough, 
 $K$ can be covered by $ r_K/\rho^{2n}$  balls of radius $\rho$,
 in such a way that every point of $K$ belongs to at most $n_K$ balls.
 \end{lemme}
 
 \bpr The lattice $\Zz^{2n}$ acts by translations on $\Cc^n$. 
 The orbit of the  ball $B^{2n}(\sqrt n)$
  under this action   covers $\Cc^n$ in such a way that every point 
  belongs to a finite number of balls. The images of this covering under
  homothetic transformations provides for every  $\rho>0$ a covering of $\Cc^n$ 
  by balls of radius $\rho$ such that every point belongs to a 
  number of balls bounded independently of $\rho$. 
  Let  $(X,\omega)$ be a K\"ahler manifold. For every point $x\in K$, 
  choose a holomorphic embedding  $\phi_x : B' \to X$, where
  $B'$ is a ball in $\Cc^n$ independent of $x$, 
  $\phi_x(0) = x$ and $\phi_x$ is everywhere contracting.
  Let $B\subset B'$ be the ball of half radius. We extract
  a finite subcovering $\phi_1(B), \cdots, \phi_k(B)$ from 
  the covering $(\phi_x(B))_{x\in K}$ of $K$. For every 
  $j \in \{1, \cdots, k\}$ and every $p\in B$, there exists an affine expanding
  map $D^j_p : \Cc^n\to \Cc^n$ that fixes $p$ and such that
  $\phi_j\circ D^j_p$ is an isometry at $p$. Then, there exists  $\rho_0>0$
  such that for every $0<\rho \leq \rho_0$ and  $1\leq j \leq k$, the restriction
  to $B$ of the covering of $\Cc^n$ by balls of radius $\rho$ satisfies the following:
  for every ball $B_p(\rho)$ of this covering centered at $p\in B$, 
  we have $D^j_p (B_p(\rho))\subset B'$, and $\phi_j\circ D^j_p(B_p(\rho))$ is 
  a ball of radius $\rho$ of $(X,\omega)$. Since $D^j_p$ is expanding, 
  $D^j_p(B_p(\rho))$ contains $B_p(\rho)$, 
  so that the union of these balls of radius $\rho$ covers $K$. 
  Moreover, the norm of $D^j_p$ is uniformly bounded  on $B$. 
  Thus, there exists a constant $h>1$ such that $D^j_p(B_p(\rho))\subset B_p(h\rho)$
  for every $p$ and $j$. From this and the construction of our coverings of $\Cc^n$, we deduce the existence of a constant $n_K>0$ independent of $p$ such that for every point $x\in K$
  and every covering of $K$ by balls of radius $\rho$ obtained in this way, 
  $x$ belongs to  at most $n_K$ balls of the covering. Finally, the existence
  of $r_K$ follows from the construction of the  covering of $\Cc^n$ we used. 
  \epr
 
\section{Proof of the theorems}
\subsection{Proof of Theorem \ref{dimension1}}
\begin{lemme}\label{chi-eta}
Let $X$ be a closed real  one-dimensional K\"ahler manifold.
There exist nonnegative constants $\eta_0$, $E_1$, $E_2$, $E_3$, $E_4$
and a family of smooth cutoff functions $\chi_\eta : X \to [0,1]$ with support
in $\RX$, $0<\eta \leq \eta_0$, such that for every $0<\eta\leq \eta_0$,
\begin{enumerate}
\item $ E_1 \eta \leq Vol (supp (\partial \del \chi_\eta))$
\item $Vol (X\setminus \chi_\eta^{-1}(1) ) < E_2 \eta$
\item $||\partial \del \chi_\eta||_{L^\infty} \leq E_3/\eta^2$
\item $dist (supp (\chi_\eta), \Rr X)\geq E_4 \eta.$
\end{enumerate}
\end{lemme}

\bpr
A neighborhood $V$ of the real locus $\Rr X$ is the union
of a finite number of annuli  isomorphic to  $A = \{ z\in \Cc \ | \ 1-\epsilon < |z| < 1+\epsilon       \}$.
For every $\eta>0$, choose $\chi_\eta$ such that $\chi_\eta(X\setminus V) = 1$
and the restriction of $\chi_\eta$ to $A$  only depends on the modulus of $z\in A$. 
That is, for every $z\in A$, $\chi_\eta(z) = \rho_\eta(|z|-1)$, where
$\rho_\eta$ is a function $ ]-\epsilon , \epsilon[ \to [0,1]$. Let 
$ \rho : \Rr \to [0,1]$ be an even function such that $\rho(x) = 1$ if
$ |x|\geq 1$ and $\rho(x) = 0$ if $|x| \leq 1/2$. For every $\eta>0$,
we set $\rho_\eta(x) = \rho(x/\eta)$. The family $\chi_\eta$, $0<\eta \leq \epsilon =\eta_0$
satisfies the required conditions.
\epr

\bpr[ of Theorem \ref{dimension1}]
For every $d\in \Nn^*$, denote by
$$ \bm_d^{\epsilon (d)} = \{ \si \in \Rhd \setminus \Rr \Delta_d \  | \ \#(\si^{-1}(0)\cap \Rr X)\geq \sqrt d \epsilon (d)       \}. $$
For every $\si \in \bm^{\epsilon (d)}_d$, denote by $Z_\si : \phi \in C^0(X) \mapsto \frac{1}{d}\int_{C_\si} \phi \in \Rr$ the associated discrete measure, where $C_\si = \si^{-1}(0)$. 
Let $(\chi_\eta)_{0<\eta \leq \eta_0}$ be a family of real cutoff functions given by Lemma \ref{chi-eta}.
By definition,  for every $0<\eta\leq \eta_0$ and every $\si \in \bm_d^{\epsilon(d)}$,
$$ \langle Z_\si, \chi_\eta \rangle = \frac{1}{d}\int_{C_\si} \chi_\eta  \leq \int_X \omega - \frac{\epsilon(d)}{\sqrt d},$$
where $\omega$ denotes the curvature of $L$. Without loss of generality, we can assume that
when $d$ is large enough, $\epsilon(d)/\sqrt d \leq \eta_0$.  We then set $ \eta_d = \epsilon (d)/(2E_2\sqrt d \int_X \omega),$
 where $E_2$ is given by Lemma \ref{chi-eta}.  From Lemma \ref{chi-eta}, we deduce that
 for every $\si \in \bm^{\epsilon (d)}_d$, 
$\langle \omega- Z_\si, \chi_{\eta_d}\rangle > \frac{\epsilon (d)}{2\sqrt d}$
and then from Poincar\'e-Lelong formula that 
$$ \frac{1}{d}|\int_X \log ||\si(x) ||^2_{\Phi_d} \partial \del \chi_{\eta_d}| \geq \frac{\pi \epsilon (d)}{\sqrt d} - 
2\pi || \frac{1}{d}\Phi^*\omega_{FS} - \omega||_{L^\infty}.$$ 
We know from Tian's asymptotic isometry theorem  \cite{Tian} that $d|| \frac{1}{d}\Phi^*\omega_{FS} - \omega||_{L^\infty}$
is bounded, so that for $d$ large enough, the right hand side is greater than $\epsilon(d)/\sqrt d$. 
For every $d$ large enough, denote by  $K_d $ the  support of $ \partial \del \chi_{\eta_d}$. 
Without loss of generality, we can assume that  $\epsilon(d)$ grows to infinity
when $d$ grows to infinity. By Lemma \ref{chi-eta}, so does $d \text{ dist}(K_d, \Rr X)^2$. 
Proposition \ref{prop} and Lemma \ref{chi-eta} then provide the result.
\epr
 \subsection{Proof of  Theorem  \ref{principal} when $a$ is a bounded function}\label{bounded}
Let $a\in \Qq^*_+$. We have to prove  the existence of  two positive constants
$C$ and $D$ such that $ \mu(\bm^a_d) \leq C e^{-Dd}.$
From Theorem \ref{laminaire} we know that the compact
 $\bz^a = \overline{\bigcup_{d\in \Nn^*}  \bz^a_d}$ introduced in \textsection \ref{paragraph2}
 is included in $Lam(X\setminus \Rr X)$, whereas from Lemma \ref{lemme},
 $\omega \notin  Lam(X\setminus \Rr X)$. As a consequence, 
there exists a finite set $(\phi_j)_{j\in \bj}$ of two-forms with compact support in $X\setminus \Rr X$ 
such that $ \forall T \in \bz^a, \exists j\in \bj, \  | \langle \omega -T,\phi_j\rangle| > 1.$
Moreover, Poincar\'e-Lelong formula writes
$$ \forall d\geq d_L, \ \forall \si \in \Rhd \setminus \{0\}, 
\frac{1}{2i\pi d }\partial \del \log || \si  ||^2_{\Phi_d} = \frac{1}{d} \phi_d^*\omega_{FS} - Z_\si,$$
where $\omega_{FS}$ denotes the Fubini-Study form of $P(\hd)^*$ defined in 
\textsection \ref{paragraph}, and $Z_\si$  the current of integration
defined in \textsection \ref{paragraph2}. From  Tian's asymptotic isometry theorem \cite{Tian}
(see also \cite{Catlin} and \cite{Zelditch}), $\frac{1}{d} \Phi_d^*\omega_{FS}$ converges
to $\omega$ as $d$ grows to $\infty$. Thus, there exists  $d_1\geq d_L$ 
such that 
$$ \forall d\geq d_1, \forall \si \in \bm^a_d, \exists j\in \bj, \ 
\left| <\frac{1}{d}\partial \del  \log ||\si||^2_{\Phi_d}, \phi_j>\right| >2\pi.       $$
From this  relation and Proposition \ref{prop} we deduce 
\beq
 \mu(\bm^a_d) &\leq &\sum_{j\in \bj} \mu \left \{
      \si \in \Rhd \setminus \{0\} \ | \ \frac{1}{d} \left| \int_X\log ||\si||_{\Phi_d} \partial \del \phi_j\right|>2\pi \right \}\\
      & \leq &\frac{ 2c_K N_d\# \bj }{\inf Vol (supp (\phi_j))} 
                  \exp \left(- \frac{\pi d }{\max_{j\in \bj} ||Ê\partial \del \phi_j ||_{L^\infty} Vol (K) }\right),
                  \eeq
         where $K \subset X \setminus \Rr X$ is a compact containing all supports
         of the $\phi_j$'s, $j\in \bj$, and $ c_K$ is given by Corollary \ref{corollaire}. 
      Hence the result. $\square$

     \subsection{Proof of Theorem \ref{principal}, general case}
     \begin{lemme}\label{lemme3}
     Under the hypotheses of Theorem \ref{principal}, let $\stackrel{\circ}{K}$ be
 a relatively compact open subset of $X \setminus \Rr X$ equipped with  
a covering by balls given by  Lemma \ref{lemme2}. Let $n_K$ be given by Lemma \ref{lemme2}
and $\lambda_K>0$, $\ba^d_K\subset \Rhd$ be given by  Corollary \ref{corollaire2}. 
Then for every $ d\geq d_L$ and  $\si \in \bm^{a(d)}_d \setminus \ba^d_K$, there is a ball $B$
of the covering such that the genus $g(C_\si\cap B)$ of $C_\si\cap B$
satisfies  
$g(C_\si\cap B) \leq \frac{n_K}{\lambda_K} a(d) A(C_\si\cap B),$ 
where $A(C_\si \cap B)$ 
denotes the area of $C_\si\cap B$.
     \end{lemme}
     Recall that the integer $d_L$ was defined in \S \ref{asymptotic}.

\bpr By definition, the genus $g(C_\si \cap B)$ is such that the Euler characteristic
 $\chi(C_\si \cap B)$ of this curve be given by the formula
 $$ \chi(C_\si \cap B)=  2b_0(C_\si \cap B) - 2g(C_\si \cap B ) - r(C_\si \cap B)$$
 where $b_0(C_\si \cap B)$ (resp. $r(C_\si \cap B)$) denotes the number of the connected components
 of $C_\si \cap B$ (resp. of $\partial (C_\si \cap B)$). In particular, for every $\si \in \bm^{a(d)}_d, \ g(C_\si \cap B) \leq g(C_\si \setminus \Rr C_\si) \leq a(d)d.$
 Denote by $(B_i)_{i\in \bi}$ the covering of $K$. Since $\si \notin \ba^d_K$,  Corollary \ref{corollaire2} implies
 that 
$ \sum_{i\in \bi} A(C_\si \cap B_i) \geq A(C_\si \cap K) \geq \lambda_Kd$. 
 Now, from Lemma \ref{lemme2} we conclude that
 $$  \sum_{i\in \bi} g(C_\si \cap B_i) \leq n_K g (C_\si \setminus \Rr C_\si) \leq n_K a(d)d.$$
 Hence the result.  
  \epr

\bpr[ of Theorem \ref{principal}]
From \textsection \ref{bounded}, we can assume that the sequence $a(d)$ grows to infinity. 
For every $d\in \Nn^*$, we set $\rho_d = a(d)^{-1/2}.$
Let $\stackrel{\circ}{K}$
be  a relatively compact  open subset of $X \setminus \Rr X$. For every
 $d$ large enough, we cover $K$ by balls of radius $\rho_d$
  as given by Lemma \ref{lemme2}.
 This cover contains at most $r_K/\rho^4_d = r_K a(d)^2$ balls. 
 From Corollary \ref{corollaire3}, there is a subset $\bb^d$ of 
 $ \Rhd$  satisfying 
 $$ \mu(\bb^d) \leq2\int_X \omega^2 c_K r_K N_d a(d)^4 \exp \left(-D_K \frac{d}{a(d)}\right), $$ 
 and such that for every $\si \in \Rhd \setminus \bb^d$ and every ball $B$ of the cover,
 $$ \lambda^1_K \frac{d}{a(d)^2} \leq A(C_\si \cap B)\leq \lambda_K^2 \frac{d}{a(d)^2}.$$
 Let $\ba^d_K \subset \Rhd$ be the set given by Corollary \ref{corollaire2}. By Lemma \ref{lemme3},
 for every $\si \in \Rhd \setminus \ba^d_K$, there is a ball $B_\si$ of our cover such that
 $$ g(C_\si \cap B_\si) \leq \frac{n_K}{\lambda_K}a(d) A(C_\si \cap B).$$ 
 Without loss of generality, we can assume that $B_\si = B^4(\rho_d)\subset \Cc^2$
 and that the area of $C_\si$ is computed with respect to the standard metric $\omega_0$ 
 of $\Cc^2$. Denote by  $\widetilde C_\si $ the image of $C_\si  \cap B$ 
under the homothetic transformation with coefficient $1/\rho_d$, so that $\widetilde C_\si \subset
 B^4(1)$ and $g(\widetilde C_\si) \leq \frac{n_K}{\lambda_K} A(\widetilde C_\si).$
 
Let  $\bt^{(1,1)}_{\pi^2} (B(1))$ be the space of positive closed currents of bidegree $(1,1)$ 
 on the unit ball $B^4(1)$ with mass $\pi^2$, and $\widetilde Z_\si \in \bt^{(1,1)}_{\pi^2} (B(1))$
 the current of integration 
 $$ \widetilde Z_\si : \phi \in \Omega^{(1,1)}_c (B^4(1)) \mapsto \widetilde Z_\si (\phi) = \frac{\pi^2}{A(\widetilde C_\si)} \int_{\widetilde C_\si} \phi.$$
 We set  $$ \bz^a =  \overline{\bigcup_{d\geq d_L}Ê\{\widetilde Z_\si, \si \in \Rhd \setminus \ba^d_K \}}  \subset \bt^{(1,1)}_{\pi^2} (B(1)). $$
 By  Theorem 1 of \cite{deThelin}, $\bz^a$ is contained in the space
 of weakly laminar currents of the unit ball $B^4(1)$. In particular, from Lemma \ref{lemme}
 we know that $\omega_0 \notin \bz^a$. Since $\overline{B^4(1)}$ 
 is compact, $\bz^a$ is compact and there exists a finite number of two-forms $(\widetilde \phi_j)_{j\in \bj}$
 with compact support in $B^4(1)$ such that 
 $$ \forall \lambda \in [\frac{\lambda^1_K}{\pi^2}, \frac{\lambda^2_K}{\pi^2}], \forall  T\in \bz^a, \exists j \in \bj, Ê|\langle \lambda T- \omega_0,\widetilde \phi_j\rangle| > 1.$$
 Applying this inequality to $T= \widetilde Z_\si$ and $\lambda = a(d)^2 A(C_\si \cap B)/\pi^2 d$, we get
 $$ \left|  \frac{a(d)^2 A(C_\si \cap B)}{ dA(\tilde C_\si)} \int_{\widetilde C_\si} \widetilde \phi_j -Ê\int_{B(1)} \omega_0 \wedge \widetilde \phi_j\right|Ê> 1.$$
 Denote by $\phi_j$ the pullback of $\widetilde \phi_j$ under
  the homothetic tranformation of coefficient $1/\rho_d$,
 so that the support of $\phi_j$ lies in $B^4(\rho_d)$. We get 
 $$ \left|  \frac{1}{d} \int_{ C_\si} \phi_j -Ê \int_{B^4(\rho_d)} \omega_0 \wedge \phi_j\right|Ê> 1/a(d),$$
 as long as $\si \notin \bb^d$. 
Finally, since by Definition \ref{ball} $ a(d) \int_{B_\si} (\omega - \omega_0) \wedge \phi_j$ 
converges to zero, we deduce   for $d$ large enough the relation
$$ \forall \si \in \Rhd \setminus (\ba^d_K\cup \bb^d), \ \exists j\in \bj, \
Ê\left| \frac{1}{d} \int_{ C_\si} \phi_j -Ê\int_X \omega \wedge \phi_j\right|Ê\geq 1/ a(d).$$
Likewise, from Tian's asymptotic isometry theorem \cite{Tian}, $a(d) \int_X (\omega - \frac{1}{d} \Phi_d^*\omega_{FS})
\wedge \phi_j$ converges to zero as $d$ grows to $\infty$. 
Applying Proposition \ref{prop}  to every ball of our cover and every $\phi_j$, $j\in \bj$,
with support in this ball, we finally obtain the existence of positive constants $C$, $D$, such that
$ \mu(\bm^{a(d)}_d) \leq CN_d a(d)^4 \exp \left(- D \frac{d}{a(d)}\right).$
\epr
\section{Final remarks}
\subsection{Average current of integration}
For every $k\geq 1$, denote by
\beq E_{\Cc P^k} : \Cc P^k & \to & \Rr\\
z & \mapsto & \int_{\Rr H^0(\Cc P^k, \bo_{\Cc P^k}(1))} \log ||Ê\si(z)||^2 d\mu(\si)
\eeq
the expectation of the random variable $ \si \mapsto \log ||\si||^2$. 
\begin{prop}\label{propE}
For every $k\geq 1$ and  $z \in \Cc P^k \setminus \Rr P^k$, 
$$ E_{\Cc P^k}(z) = \log (\frac{k+1}{4}) + \int^\infty_0 e^{-\rho} \log \rho d\rho + \log (1 + \sqrt {1 - || \tau||^2(z)}),$$
where $\tau $ is the section introduced in Proposition \ref{tau}.
\end{prop}
This result is very close to Lemma 2.5 of \cite{Macdonald}.

\bpr As in the proof of Proposition \ref{tau}
and using the notations of Remark \ref{invariance}, we get
for every  $0 < r\leq 1$:
\beq E_{\Cc P^k}(z_r) &=& \int_{\Rr^2} \log \left((k+1) \frac{|a_0 + ira_1|^2}{1+r^2}\right) \frac{e^{-|a|^2}}{\pi} da_0 da_1\\
&= &\log \left(\frac{k+1}{1+r^2}\right) +  \int^\infty_0  e^{-\rho} \log \rho d\rho  + 
\frac{1}{2\pi} \int_0^{2\pi} \log |\cos \theta +ir\sin \theta|^2 d\theta\\
& = & \log \left(\frac{k+1}{4(1+r^2)}\right) +  \int^\infty_0  e^{-\rho} \log \rho d\rho  + 
\frac{1}{2\pi} \int_0^{2\pi} \log |e^{2i\theta}(1+r) + 1-r |^2 d\theta.
\eeq
From Jensen formula, as soon as $r>0$, 
$$ \frac{1}{2\pi} \int_0^{2\pi} \log |e^{2i\theta}(1+r) + 1-r |^2 d\theta= \log |1-r|^2 + \log \left|\frac{1+r}{1-r}\right|^2 = \log (1+r)^2. $$
From this we deduce 
 \beq E_{\Cc P^k}(z_r) &=&\log (\frac{k+1}{4}) +  \int^\infty_0  e^{-\rho}\log \rho d\rho  +\log \left(\frac{(1+r)^2}{1+r^2}\right)\\
&= &  \log (\frac{k+1}{4}) +  \int^\infty_0 e^{-\rho} \log \rho  d\rho  + \log \left(1 + \sqrt { 1 - || \tau||^2(z)}\right),
\eeq
since
$||\tau||(z_r) = |\tau(z_r)|/|z_r|^2 = (1-r^2)/(1+ r^2).$
The result follows from the invariance of $E_{\Cc P^k}$ and $||\tau||$ under the action
of  $PO_{k+1}(\Rr)$, see Remark \ref{invariance}.
\epr

\begin{coro}\label{ddbar}
For every $k\geq 1$ and every real line $D$ in $\Cc P^k$, the restriction of the current  $\frac{1}{2i\pi} \partial \del E_{\Cc P^k}$
to $D\setminus \Rr D$ coincides with the Fubini-Study form, while its restriction to the quadric $\{\tau = 0\}$ vanishes.
\end{coro}

\bpr 
Proposition \ref{propE} implies that the restriction of $ E_{\Cc P^k}$ to the quadric $\{\tau = 0\}$ 
is constant, so that the current $\partial \del E_{\Cc P^k}$ vanishes on this quadric. 
In the same way, Proposition \ref{propE}
 implies that the restriction of $\frac{1}{2i\pi} \partial \del E_{\Cc P^k}$
to $D$ does not depend on $k$. Thus,  we may assume $k= 1$ and $D = \Cc P^1$. 
Now,  every $\si \in \Rr H^0(\Cc P^1, \bo_{\Cc P^1} (1))$
 does not vanish on $\Cc P^1 \setminus \Rr P^1$, so that by definition 
\beq \forall z \in \Cc P^1 \setminus \Rr P^1, Ê\frac{1}{2i\pi} \partial \del E_{\Cc P^k} (z) & = & \int_{ \Rr H^0(\Cc P^1, \bo_{\Cc P^1} (1))}
\frac{1}{2i\pi}\partial \del \log ||\si(z)||^2 d\mu (\si)\\
& = & \int_{ \Rr H^0(\Cc P^1, \bo_{\Cc P^1} (1))} \omega_{FS}(z)  d\mu(\si) \\
 & = & \omega_{FS}(z).
 \eeq
  \epr
 Now, let $L$ be a real Hermitian  line bundle with positive curvature on a smooth real K\"ahler
 manifold $X$ of dimension $n\geq 1$. For every $d\in \Nn^*$ and
 every $(2n-2)$-form $\phi \in \Omega^{2n-2}(X)$, we denote by
 $$ Z^\phi : \si \in \Rhd \setminus \{0\} \mapsto Z^\phi_\si = \frac{1}{d} \int_{C_\si} \phi \in \Rr $$
 the associated random variable, where the space $\Rhd$ is equiped with the $L^2$ Gaussian
probability measure   $\mu$. We write $$ E_d(Z^\phi) = \int_{\Rhd} Z^\phi_\si d\mu(\si)$$
 for the expectation of this random variable, and $E_d(Z) : \phi \in \Omega^2(X) \mapsto E_d(Z^\phi) \in \Rr$ for the
 associated 
 closed positive (1,1)-current. 
 \begin{prop}\label{moyenne}
 Let $L$ be a  real Hermitian  line bundle with positive curvature  on a  smooth closed real K\"ahler manifold
  $X$. Then, for every $d\geq d_L$,  
 $$ E_d(Z) = \frac{1}{d} \Phi_d^* \omega_{FS} - \frac{1}{2i\pi d} \Phi_d ^*\partial \del E_{P(\hd^*)}. $$
 Moreover, the restriction of this current to the complement of the real locus  
 converges to $\omega$ as $d$ grows to  infinity.
 \end{prop}
 Recall that the embedding $\Phi_d : X \to P(\hd^*)$, $d\geq d_L$, and the 
 Fubini-Study form $\omega_{FS}$ of the projective space $P(\hd^*)$
were introduced in \textsection \ref{asymptotic}. 
 
 \bpr Poincar\'e-Lelong formula provides for every $\si \in \Rhd\setminus \{0\}$ the relation
 $$ \frac{1}{2i\pi d}\partial \del \log||\si||^2_{\Phi_d} = \frac{1}{d}\Phi_d^* \omega_{FS} - Z_\si.$$
 The first part of Proposition \ref{moyenne} is obtained by integration of this relation 
 on $\Rhd$. 
 Tian's asymptotic isometry theorem \cite{Tian} implies that $\frac{1}{d}\Phi_d^* \omega_{FS}$
 converges to the curvature $\omega$ of $L$. Proposition \ref{sup} combined
 with Proposition \ref{propE} imply that
 $  \frac{1}{2i\pi d} \Phi_d^*\partial \del E_{P(\hd^*)}$ converges to zero faster than every
 polynomial function in $d$, and even exponentially fast in the cases covered by Remark \ref{exp} (compare
 with \cite{Macdonald}). Hence the result.
  \epr
Note  that when the chosen probability space is the whole complex space $H^0(X, L^d)$, 
 the expectation $\int_{\hd} \log ||\si(z)||^2 d\mu(\si)$ is a function of $zÊ\in \Cc P^k$ invariant
 under the whole $PU_{k+1}(\Cc)$, thus is constant. Hence, $E(Z^\phi_\Cc) = \frac{1}{d}\Phi_d^* \omega_{FS}$, 
 see \cite{Shi-Zel}. Moreover, Shiffman and Zelditch proved in \cite{Shi-Zel3} that the law of $Z^\phi_\Cc$ 
 converges to a normal law as $d$ grows to  infinity, a result that was already obtained in dimension
 one in \cite{Sodin-Tsirelson}. It would be here of interest to understand in more details the convergence of the law of $Z^\phi$. 
 
 \subsection{Existence of real maximal  curves}
 An algebraic curve $C$ of genus $g(C)$ is said to be maximal when the number of  components 
 of its real locus coincides with $g(C)+1$, the maximum allowed by  Harnack-Klein inequality \cite{Harnack},  \cite{Klein}. Our Theorem \ref{principal} proves, in particular, that if $L$ is a real ample Hermitian line bundle over a
 real projective  surface $X$, the measure of the set of real maximal curves 
 linearly equivalent to $L^d$ exponentially decreases as  $d$ increases. When $X = \Cc P^2$, Harnack \cite{Harnack}  proved that such maximal curves exist in any degree.  The study of these curves plays  a central r\^ole in real algebraic geometry, at least since Hilbert included it in his 16th problem. Nevertheless, such curves do not always exist. 
 For instance, if $X$ is the product of two non maximal real curves, then for every ample real line bundle $L$ 
 over $X$ 
 and every $d\in \Nn^*$, the linear system $\Rhd$ contains no maximal curve.

 However, every real closed symplectic manifold $(X,\omega, c_X)$ of dimension four
 with rational form $\omega$ supports, when $d$ is large enough, symplectic real surfaces
 Poincar\'e dual to $d\omega$ and whose real locus contains at least $\epsilon d$ 
 components, where $\epsilon$ depends on the manifold $(X,\omega,c_X)$, see \cite{Gayet}.
  Applying Harnack's method to these curves, we see that there always exist even
  symplectic real surfaces whose real locus contains at least $\epsilon' d^2$ connected components.

  The following questions then arise. For every ample real line bundle $L$ on a projective real surface 
  $X$ and every $d\in \Nn^*$, denote by
  $$ m(L^d) = \sup_{\si \in \Rhd\setminus \Rr \Delta_d} b_0(\Rr C_\si)$$
  the maximal number of connected components that  a smooth real  divisor linearly equivalent to $L^d$
  may contain. 
 Then, denote by $ \epsilon(X,L) = \lim \sup_{d\to \infty} \frac{1}{d^2} m(L^d),$
 so that $0\leq \epsilon (X,L) \leq \frac{1}{2} L.L$
 by Harnack-Klein inequality and the adjunction formula. 
Is  this quantity $\epsilon (X,L)$  bounded from below by a non negative constant independant of $(X,L)$? 
Does there exist a pair $(X,L)$ with $\Rr X \not= \emptyset$ such that $\epsilon(X,L) < \frac{1}{2}L^2$? If not, what
about the quantity $ \lim \sup_{d\to \infty} \frac{1}{d}(m(L^d) - \frac{1}{2}d^2 L^2)$?

  The same questions hold within the realm of four-dimensional real symplectic  manifolds. 
  Recall that the real symplectic surfaces built in \cite{Gayet} are 
  obtained via Donaldson's method \cite{Donaldson}, so that their current   of integration converges to $\omega$
  as $d$ grows to infinity. Theorem \ref{laminaire} provides an obstruction to get  real maximal curves using this method (Donaldson's quantitative transversality gives
  another one, as observed in \cite{Gayet}). 
  This phenomenon was in fact the starting point of our work.

  This  work raises several questions.
   It is known \cite{EK} that the expectation of the number of real roots of a real polynomial in one variable 
  is $\sqrt n$.  What is the expected value  of $b_0(\Rr C_\si) $ in dimension two? 
   How to improve Theorem \ref{principal} to get decays till this expectation, as in Theorem \ref{dimension1}? 
   What happens for values below  this expectation?
   Note that for spherical harmonics on the two-dimensional sphere, 
   such kinds of results have been obtained in \cite{Nazarov-Sodin}.    
   More generally, what is the  asymptotic  law of the random variable $b_0$?  
   What happens in higher dimensions?


\noindent
\textsc{Universit\'e de Lyon ; CNRS \\
Universit\'e Lyon 1 ; Institut Camille Jordan}


\end{document}